\documentclass[12pt]{article}

\usepackage{amssymb}

\usepackage{amsfonts}
\usepackage{latexsym}
\usepackage{amsmath,amsthm}
\newtheorem{theorem}{Theorem}%[section]
\newtheorem{proposition}{Proposition}%[section]
\newtheorem{lemma}{Lemma}%[section]
\newtheorem{definition}{Definition}%[section]
\newtheorem{example}{Example}%[section]
\newtheorem*{KKT}{Karush-Kuhn-Tucker's Theorem}{\bf}{\it}
\newcommand{\R}{\mathbb R}
\newcommand{\RR}{\overline{\mathbb R}}
\newcommand{\E}{\mathbb R^s}

\newcommand{\pr}{\prime}

\newcommand{\f}{f : X\to\R^n}

\newcommand{\scalpr}[2]{\langle#1,#2\rangle}
\newcommand{\sld}[1]{#1^{\pr\pr}_{-}} 
\newcommand{\norm}[1]{\Vert#1\Vert}
\title{Second-order optimality conditions and Lagrange multiplier characterizations of the solution set in quasiconvex programming}
\author{
Vsevolod I. Ivanov
\\ \small Technical University of Varna,
 \small Department of Mathematics,
 \small 9010 Varna, Bulgaria
}
\date{\today}

\begin{document}
\maketitle

\begin{abstract}
Second-order optimality conditions for vector nonlinear programming problems with inequality constraints are studied  in this paper. We introduce a new second-order constraint qualification, which includes Mangasarian-Fromovitz constraint qualification as a particular case. We obtain necessary and sufficient conditions for weak efficiency of problems with a second-order pseudoconvex vector  objective function and quasiconvex constraints. We also  derive Lagrange multiplier characterizations of the solution set of a scalar problem with a second-order pseudoconvex objective function and quasiconvex inequality constraints, provided that one of the solutions and the Lagrange multipliers in the Karush-Kuhn-Tucker conditions are known. At last, we introduce a notion of a second-order KKT-pseudoconvex problem with inequality constraints. We derive sufficient and also necessary conditions for efficiency of second-order KKT-pseudoconvex problems.   Three examples are presented.


\end{abstract}

\section{Introduction}
\setcounter{equation}{0}
\label{s1}
Second-order optimality conditions of Karush-Kuhn-Tuc\-ker type  play essential role in development of the contemporary vector  optimization. 
In the necessary conditions, the authors suppose some constraint qualifications to ensure the vector of the Lagrange multipliers in front of the components of the objective function to be different from zero. Several second-order constraint qualifications (in short, SOCQ) have been introduced for this aim. 

In the present paper, we deal with optimality criteria of Karush-Kuhn-Tucker type for the nonlinear programming problem with inequality constraints:

\bigskip
%\noindent
Minimize$\quad f(x)\quad$subject to $\quad g(x)\leqq0$,\hfill ${\rm (VP)}$

\bigskip
\noindent
where $f:X\to\R^n$ and $g:X\to\R^m$ are given vector real-valued functions, defined on some open set $X$ in the finite-dimensional space $\R^s$.
All results given here are obtained for nonsmooth problems in terms of the standard second-order directional derivative.

In Section \ref{s3}, we derive second-order necessary conditions for a weak efficiency for problems with C$^1$ (i.e. continuously
differentiable) data. We introduce a new second-order constraint qualification. It includes a second-order derivative and extends the Mangasarian-Fromovitz constraint qualification \cite{man67}. We call it second-order  Mangasarian-Fromovitz constraint qualification (in short SOMFCQ). If the constraints satisfy the Mangasarian-Fromovitz CQ, then they satisfy the second-order one. We show that the necessary conditions without any constraint qualification are sufficient for weak global efficiency.
In our necesary and sufficient conditions, we suppose that all components of the objective function are
second-order pseudoconvex, a notion recently introduced by Ginchev and Ivanov \cite{gi05}, and the constraint functions are quasiconvex. 

In particular, we extend the result \cite[Theorem 1]{jmaa1} to the multiobjective case. We also improve several classical sufficient optimality conditions concerning first-order case (see for example, \cite[Theorem 4.2.11]{baz79}).

In 1988, Mangasarian \cite{man88} obtained characterizations of solution sets of  convex programs in terms of a known solution. Later appeared other characterizations of the solution set of the scalar nonlinear programming problem, which concern convex, pseudolinear,  pseudoconvex, invex and other types of problems. Jeyakumar, Lee and Dinh obtained Lagrange multiplier characterizations of a convex  problem with inequality constraints and more general ones, provided that the Karush-Kuhn-Tucker multipliers are also known. Recently Suzuki and Kuroiwa \cite{suz15} obtained characterizations of the solution set of quasiconvex set-constrained problem in terms of the Greensberg-Pierscalla subdifferential. Ivanov \cite{jota4} derived characterizations of the solution sets of differentiable  quasiconvex set-constrained and inequality-constrained problems.

In Section \ref{s5}, we derive characterizations of the solution set of the scalar problem with inequality constraints in terms of a known solution $\bar x$ and known Karush-Kuhn-Tucker multipliers.
%, that is the case in problem (VP), when $n=1$ and $X$ is a convex set, not necessarily open.
We suppose in them that the objective function is second-order pseudoconvex, the inequality constraints are quasiconvex, and the set-constraint is convex, not necessarily open.

In Section \ref{s4}, we introduce a notion of a pseudoconvex vector function.  We obtain necessary and sufficient conditions for global efficiency in pseudoconvex vector problems.  In particular, we extend \cite[Theorem 3.1 a]{jogo1} to the vector case.

The preprint \cite{arxiv2013} contains an initial version of this work. Especially, the sufficient conditions in Theorems \ref{th1} and \ref{th2} appeared there.

\section{Preliminaries}
\label{s2}

We begin with some preliminary definitions.
\smallskip
Denote by $\R$ the set of reals, by $\R^n_+$ the orthant in $\R^n$ with non-negative components, and let $\RR=\R\cup\{-\infty\}\cup\{+\infty\}$.
We suppose the following rule for multiplication with infinities: $0.(\pm\infty)=0$.

Let the function $\f$ with an open domain $X\subset\R^n$ be differentiable at
the point $x\in X$. Then, the second-order directional derivative
$f^{\pr\pr}(x,u)$ of $f$ at the point $x\in X$ in direction $u\in\E$
is defined as element of $\RR$ by
\begin{displaymath}
f^{\pr\pr}(x,u)=\lim_{t\to +0}\;2 t^{-2}[f(x+tu)-f(x)-t\scalpr{\nabla f(x)}{u}].
\end{displaymath}
The function $f$ is called second-order directionally differentiable on $X$,
iff the derivative $f^{\pr\pr}(x,u)$ exists for each point $x\in X$ and any direction $u\in\R^s$.

The above derivative exists for large class of nonsmooth functions, including $l_1$-functions, max-functions, and exact penalty functions \cite{btz82}. If the function belongs to the class C$^{1,1}$, then it is always finite.
%The problem (P) is said to be Fr\'echet differentiable if $f$ and $g$ are Fr\'echet differentiable.

\smallskip
Recall that a scalar function $f:X\to\R$ is said to be (semistrictly) quasiconvex at the point $x\in X$
(with respect to $X\subset\R^s$) \cite{man69}, iff the conditions 
\[
\begin{array}{c}
y\in X,\; f(y)\leqq f(x),\; t\in[0,1],\; (1-t)x+ty\in X\;\textrm{ imply }\; f((1-t)x+ty)\leqq f(x) \\
(y\in X$, $f(y)<f(x)$, $t\in (0,1)$,  $(1-t)x+ty\in X\;\textrm{ imply }\; f((1-t)x+ty)<f(x)).
\end{array}
\]
When the set $X$ is convex the function $f$ is called  quasiconvex or semistrictly quasiconvex on $X$, iff these implications are satisfied for every $x\in X$.
%when for all $x,\;y\in X$ and $t\in [0,\,1]$ it holds
%$f((1-t)x+ty)\leqq\max\left(f(x),\,f(y)\right). $
%If the set $X$ is convex, then the function $f$ is called semistrictly quasiconvex on $X$
%when for all $x,\;y\in X$ and $t\in (0,\,1)$ it holds
%$f((1-t)x+ty)<\max\left(f(x),\,f(y)\right). $

Throughout this paper, we use the following notations comparing the vectors $x$ and $y$ with components $x_i$ and $y_i$ in finite-dimensional spaces:
\[
\begin{array}{l}
x<y,\quad\textrm{if}\quad x_i<y_i\quad\textrm{for all indexes }i; \\
%\]
%\[
x\leqq y,\quad\textrm{if}\quad x_i\leqq y_i\quad\textrm{for all indexes }i; \\
%\]
%\[
x\le y,\quad\textrm{if}\quad x_i\leqq y_i\quad\textrm{for all indexes }\quad i\quad\textrm{with at least one being strict.}
\end{array}
\]
%\textit{}

%\begin{definition}
A feasible point $\bar x\in S$ is called a (weak) local Pareto minimizer, or (weakly) efficient, iff there exists a neighborhood $U\ni\bar x$ such that there is no $x\in U\cap S$ with  $f(x)\le f(\bar x)$ ($f(x)<f(\bar x)$). 
The point $\bar x\in S$ is called a (weak) global Pareto minimizer, iff there does not exist $x\in S$ with $f(x)\le f(\bar x)$ ($f(x)<f(\bar x)$).
%\end{definition}

It is obvious that every efficient point is weakly efficient. The converse is not satisfied (see, for example, the book \cite{saw85}).

Consider the problem {\rm{(VP)}}. Denote
\[
S:=\{x\in X\mid g_j(x)\le 0,\; j=1,2,\ldots,m\}.
\]
For every feasible point $x\in S$  let $A(x)$ be the set of active constraints
\[
A(x):=\{j\in\{1,2,\ldots,m\}\mid g_j(x)=0\}.
\]
A direction $d$ is called critical at the point $x\in S$, iff
\[
\scalpr{\nabla f_i(x)}{d}\leqq 0\quad\textrm{for all}\quad i\in\{1,2,\ldots,n,\}\quad\textrm{and}\quad\scalpr{\nabla g_j(x)}{d}\leqq 0\quad\textrm{for all}\quad j\in A(x).
\]

Consider the following scalar problem:

\bigskip
Minimize$\quad h(x)\quad$subject to$\quad g_i(x)\leqq0$, $\quad i=1,2,...m$,\hfill ${\rm (SP)}$

\bigskip
\noindent
where $h:X\to\R$, $g_i:X\to\R$, $i=1,2,...,m$ are given real-valued functions, defined on some open set $X\subset\R^s$. 

The following theorem, which contains second-order necessary conditions for the scalar problem (SP), is due to Ginchev and Ivanov \cite{jmaa1}:

\begin{lemma}[]\label{NKT}
Suppose that $X$ is an open set in the space $\R^s$, and the functions
$h$, $g_i$, $(i=1,2,...,m)$ are defined on $X$.
Let the feasible point $\bar x$ be a local minimizer of the problem {\rm{(SP)}},
and the functions $h$, $g_i$ belong to the class C$^1$, $i\in A(\bar x)\}$.
Suppose that they are second-order directionally differentiable at $\bar x$ in every critical direction $d\in\R^n$,
and the functions $g_i$ $(i\notin A(\bar x))$ are continuous at $\bar x$.
%the functions $h_i$ $(i\in\{0\}\cup A(\bar x))$ are continuously differentiable,
%and the functions $h_i$ $(i\in I_0(\bar x,d))$ are second-order directionally differentiable .
Then corresponding to any critical direction $d$ there exist non-negative multipliers $\lambda,\mu_1,...,\mu_m$, not all zero, such that
\begin{gather}
\mu_i g_i(\bar x)=0,\;i=1,2,...,m,\quad\nabla L(\bar x)=0, \notag \\
%\]
%\[
\lambda\nabla h(\bar x)=0,\quad\mu_i\scalpr{\nabla g_i(\bar x)}{d}=0,\quad i\in A(\bar x), \notag \\
%\]
%\begin{eqnarray}
%\tag{4}
L^{\pr\pr}(\bar x,d)=h^{\pr\pr}(\bar x,d)+\sum_{i\in A(\bar x)}\mu_i g_i^{\pr\pr}(\bar x,d)\geqq 0. \notag
\end{gather}
\end{lemma}

%For every critical direction $d$ at the feasible point $x$ denote by $I(x,d)$ and $J(x,d)$ the following index sets:
%\[
%\centerline{$I(x,d):=\{i\in\{1,2,\ldots,n\}\mid\nabla f_i(x)d=0\},\; J(x,d):=\{j\in A(x)\mid\nabla g_j(x)d=0\}.$}
%\]

The following result is known and it could be found, for instance, in the book \cite{man69} (see Theorem 9.1.4]).
\begin{lemma}\label{lema3}
Let $X$ be an open set in $\R^s$, and let $f$ be a real scalar function, defined on $X$, which is both differentiable and quasiconvex at the point $x\in X$. Then, the following implication holds:
\[
y\in X,\,f(y)\leqq f(x)\quad\Longrightarrow\quad\scalpr{\nabla f(x)}{y-x}\leqq0.
\]
\end{lemma}

\smallskip
Let the scalar function $f: X\to\R$ with an open domain $X\subset\R^s$ be differentiable at the point $x\in X$. Then, $f$ is said to be pseudoconvex at $x\in X$, iff $y\in X$ and $f(y)<f(x)$ imply $\scalpr{\nabla f(x)}{y-x}<0$. If $f$ is differentiable on $X$, then it is called
pseudoconvex on $X$, when $f$ is pseudoconvex at each $x\in X$.

\smallskip
The following definition is due to Ginchev and Ivanov \cite{gi05}. 
\begin{definition}
Consider a function $\f$ with an open domain $X$, which is differentiable at $x\in X$ and second-order directionally differentiable at $x\in X$ in every direction $y-x$ such that $y\in X$, $f(y)<f(x)$, $\scalpr{\nabla f(x)}{y-x}=0$. Then, $f$ is called  second-order pseudoconvex  at $x\in X$, iff  for all $y\in X$ the following implications hold:
\[\begin{array}{c}
f(y)<f(x)\quad\mbox{implies}\quad\scalpr{\nabla f(x)}{y-x}\leqq0; \\
%\]
%\end{equation}
%\[
f(y)<f(x),\;\scalpr{\nabla f(x)}{y-x}=0\quad\mbox{imply}\quad f^{\pr\pr}(x,y-x)<0.
\end{array}\]

Suppose that  $f$ is differentiable on $X$ and second-order directionally differentiable at every $x\in X$ in each direction
$y-x$ such that $y\in X$, $f(y)<f(x)$, $\scalpr{\nabla f(x)}{y-x}=0$. Then, $f$ is second-order pseudoconvex on  $X$, iff it is
second-order pseudoconvex at every $x\in X$.
\end{definition}
 It follows from this definition that every differentiable pseudoconvex function is second-order pseudoconvex. The converse does not hold.

The following result is a particular case of Theorem 4 in \cite{gi05}:

\begin{lemma}\label{lema2} 
Every radially lower semicontinuous second-order pseudoconvex function, which is defined on some convex set $X\subset\R^s$, is semistrictly quasiconvex on $X$, and moreover, it is quasiconvex on $X$.
\end{lemma}

Consider the following sets:

\begin{equation}\label{1}
C_k=\{x\in X\mid g(x)\leqq 0,\; f_i(x)\le f_i(\bar x), i\ne k\}.
\end{equation}

The following result is due to Luc and Schaible \cite{luc97}:

\begin{lemma}\label{lema1}
Let the functions $f_i$ $(i=1,2,\ldots,n)$ be quasiconvex and semistrictly quasiconvex. Then, the point $\bar x$ is a weak Pareto minimizer, if and only if there exists an index $k\in\{1,2,\ldots,n\}$ such that $\bar x$ minimizes some component $f_k$ of the vector function $f$ over the constraint set $C_k$, defined by equation {\rm (\ref{1})}.
\end{lemma}

It is said that Mangasarian-Fromovitz constraint qualification \cite{man67} holds at the point $\bar x$, iff there exists a vector $u$ such that 
\[
\scalpr{\nabla g_j(\bar x)}{u}<0,\quad\forall\; j\in A(\bar x).
\]

\section{Necessary and sufficient conditions for a weak global minimum}
\label{s3}
In this section, we derive necessary and sufficient optimality conditions for weak  efficiency in the problem {\rm{(VP)}}.

We introduce the following more general  constraint qualification:

\begin{definition}\label{df1}
The following condition is an extension of Mangasarian-Fromovitz constraint qualification: there exists a direction $d\in\R^n$ and a set of indexes $K$, $K\subset A(\bar x)$ such that
\begin{equation}\label{2}
\begin{array}{l}
\scalpr{\nabla g_j(\bar x)}{d}=0,\quad g^{\pr\pr}_j(\bar x,d)<0, \quad j\in K \\
\scalpr{\nabla g_j(\bar x)}{d}<0,\quad j\in A(\bar x)\setminus K.
\end{array}
\end{equation}
\end{definition}

The constraint qualification from Definition \ref{df1} 
can be easy verified. We call it second-order Mangasarian-Fro\-mo\-vitz constraint qualification (in short, SOMFCQ), because second-order derivatives appear in SOMFCQ.

\begin{theorem}\label{th1}
Let the point $\bar x$ be feasible for the problem {\rm (P)}. Suppose that the set $X$ is convex, the vector function $f$ and the scalar functions $g_j$, $j\in A(\bar x)$ are continuously differentiable, the functions $g_j$, $j\notin A(\bar x)$ are continuous at $\bar x$, $f$ and $g$ are second-order directionally differentiable at $\bar x$ in every critical direction $d$, all components of $f$ are 
second-order pseudoconvex, all components of $g$ are quasiconvex. 

Let $\bar x$ be a weak Pareto minimizer. Then, for every direction $d\in\R^s$ there exist Lagrange multipliers $\lambda\in\R^n_+$, 
$\mu\in\R^m_+$ with $(\lambda,\mu)\ne (0,0)$ such that
\[
\mu_i g_i(\bar x)=0,\;i=1,2,...,m,
\]
and the following conditions are satisfied:
\begin{equation}\label{9}
%\mu_i g_i(\bar x)=0,\;i=1,2,...,m,\quad
\scalpr{\nabla L(\bar x)}{d}>0,
\end{equation}
if $d$ is not critical, or
\begin{equation}\label{10}
\nabla L(\bar x)=0,\quad L^{\pr\pr}(\bar x,d)=\sum_{i=1}^n\,\lambda_i f_i^{\pr\pr}(\bar x,d)+\sum_{j\in A(\bar x)}\mu_j g_j^{\pr\pr}(\bar x,d)\geqq 0,
\end{equation}
if $d$ is critical. If we suppose additionally that SOMFCQ holds, then 
%for every direction satisfying SOMFCQ 
there exists a direction $d$ and vector multipliers $\lambda$, $\mu$ with non-negative components such that (\ref{9}) and (\ref{10}) hold with $\lambda\ne 0$. 
%Really, for every direction satisfying SOMFCQ there exists a multiplier $\lambda$, $\lambda\ne 0$.
Here, $L$ is the Lagrange function $L=\sum_{i=1}^n\,\lambda_i f_i+\sum_{j=1}^m\,\mu_j g_j$.

Conversely, suppose that for every critical direction $d\in\R^s$ there exist Lagrange multipliers 
\[
\lambda=(\lambda_1,\lambda_2,\ldots ,\lambda_n)\in\R^n_+,\; \mu=(\mu_1,\mu_2,\ldots,\mu_m)\in\R^m_+\;\textrm{ with }\; \lambda\ne 0
\]
 such that 
conditions (\ref{10}) are satisfied. Then, $\bar x$ is a weak global Pareto minimizer.
\end{theorem}
%conditions (\ref{9}) and (\ref{10}) are satisfied.
%\begin{equation}\label{2}
%\mu_j g_j(\bar x)=0,\;j=1,2,...,m,\quad\nabla L(\bar x)=0,
%\end{equation}
%\begin{equation}\label{3} 
%L^{\pr\pr}(\bar x,d)=\sum_{i=1}^n\lambda_i f_i^{\pr\pr}(\bar x,d)+ \sum_{j\in A(\bar x)}\mu_j f_j^{\pr\pr}(\bar x,d) \geqq 0.
%\end{equation}

%{\bf Proof.}
\begin{proof}
Let $\bar x$ be a weak minimizer. 
We prove that conditions (\ref{9}) and (\ref{10}) hold. It follows from Lemma \ref{lema2} that all components of $f$ are semistrictly quasiconvex and quasiconvex. Then, it follows from Lemma \ref{lema1} that there exists an index $k\in\{1,2,\ldots,m\}$ such that $\bar x$ minimizes the function $f_k$ over the set $C_k$.  %According to Lemma \ref{NKT} conditions (\ref{9}) and (\ref{10}) are satisfied.

We prove that the second-order KKT conditions are satisfied.
Let $\scalpr{\nabla f_i(\bar x)}{d}>0$ for some index $i\in\{1,2,\dots,n\}$. Then the condition (\ref{9}) is fulfilled with 
\[
\mu_1=\mu_2=\cdots=\mu_m=0,\; \lambda_k=0,\;\textrm{ when }\;k\ne i,\; \lambda_i=1.
\]
 Consider the case, when $\scalpr{\nabla f_i(\bar x)}{d}\le 0$ for all indexes $i=1,2,\dots,n$. It is possible that for some index $j\in A(\bar x)$ is satisfied the inequality $\scalpr{\nabla g_j(\bar x)}{d}>0$. In this case, we could take $\mu_j$ to be a sufficiently large positive number and $\mu_k=0$,  when $k\ne j$, $\lambda_i=1$ for all $i=1,2,\dots,n$. This choice will ensure condition (\ref{9}). Otherwise, we have $\scalpr{\nabla g_j(\bar x)}{d}\le 0$, $j\in A(\bar x)$, which implies that $d$ is a critical direction. 
Then, the Fritz-John type second-order conditions
\begin{equation}\label{7}
\nabla L(\bar x)=0,\quad\sum_{i=1}^n\,\lambda_i f_i^{\pr\pr}(\bar x,d)+\sum_{j\in A(\bar x)}\mu_j g_j^{\pr\pr}(\bar x,d)\ge 0,\quad (\lambda,\mu)\ne (0,0)
\end{equation}
 follow directly from Lemmas \ref{lema1} and \ref{NKT}, taking into account that the constraints $f_j(x)-f_j(\bar x)\leqq 0$ are active at $\bar x$.

Let SOMFCQ hold. We prove that there exists a direction $d$ and  multipliers $\lambda$, $\mu$ such that $\lambda=\lambda(d)\ne 0$. Suppose the contrary that $\lambda=0$ for every direction $d$. Let $\bar d$ be the direction, which satisfies SOMFCQ. It is impossible that $\scalpr{\nabla f_i(\bar x)}{\bar d}>0$ for some index $i\in\{1,2,\dots,n\}$, because the choice $\mu_1=\mu_2=\cdots=\mu_m=0$, $\lambda_k=0$ when $k\ne i$, $\lambda_i=1$ will ensure $\scalpr{\nabla L(\bar x)}{\bar d}>0$ and $\lambda\ne 0$, which is impossible according to our assumption. Therefore, $\scalpr{\nabla f_i(\bar x)}{\bar d}\le 0$ for all indexes $i=1,2,\dots,n$. It follows from here that $\bar d$ is critical, because it satisfies SOMFCQ.
By Fritz John conditions (\ref{7}) and $\lambda=0$, we conclude that 
%for every critical direction 
there exists a multiplier 
$\mu=(\mu_1,\mu_2,\ldots,\mu_m)\in\R^m_+$ with $\mu\ne 0$ such that 
\begin{equation}\label{3}
\mu_j\scalpr{\nabla g_j(\bar x)}{\bar d}=0,\; j\in A(\bar x),
\end{equation}
\begin{equation}\label{5}
\sum_{j\in A(\bar x)}\mu_j g_j^{\pr\pr}(\bar x,\bar d) \geqq 0.
\end{equation}
Denote by $K_1(d)$ the index set, which depend on $d$, such that
\[
K_1(d):=\{j\in A(\bar x)\mid\scalpr{\nabla g_j(\bar x)}{d}=0\}.
\]
According to the SOMFCQ 
%second-order Man\-ga\-sa\-rian-Fromovitz CQ 
the direction $\bar d$ and the set of indexes $K(\bar d)$ satisfy the conditions
\begin{equation}\label{6}
\begin{array}{l}
\scalpr{\nabla g_j(\bar x)}{\bar d}=0,\quad g^{\pr\pr}_j(\bar x,\bar d)<0, \quad j\in K(\bar d) \\
\scalpr{\nabla g_j(\bar x)}{\bar d}<0,\quad j\in A(\bar x)\setminus K(\bar d).
\end{array}
\end{equation}
%It is easy to see that $\bar d$ is critical, and 
By SOMFCQ,
%second-order Man\-ga\-sa\-rian-Fromovitz CQ, 
we could suppose without loss of generality that $K_1(\bar d)=K(\bar d)$. By the equations (\ref{3}), we have $\mu_j(\bar d)=0$ for all $j\in A(\bar x)\setminus K(\bar d)$. It follows from $\mu(\bar d)\ne 0$ that $K(\bar d)\ne\emptyset$, and for every $j\in K(\bar d)$ such that $\mu_j(\bar d)\ne 0$, we have
$\scalpr{\nabla g_j(\bar x)}{\bar d}=0$, $g_j^{\pr\pr}(\bar x,\bar d)<0$. This is a contradiction to the condition
\[
\sum_{j\in A(\bar x)}\mu_j(\bar d) g_j^{\pr\pr}(\bar x,\bar d) \geqq 0,
\]
which follows from (\ref{5}).

Conversely, suppose that for every critical direction $d$ there exist $\lambda$ and $\mu$, which satisfy conditions (\ref{10}). We prove that $\bar x$ is a weak global Pareto minimum. Assume the contrary that there exists $x\in S$ with $f(x)<f(\bar x)$. We prove that $x-\bar x$ is a critical direction.
By second-order pseudoconvexity, $\scalpr{\nabla f_i(\bar x)}{x-\bar x}\leqq0$ for all $i=1,2,\dots,n$.
Due to quasiconvexity and $g_j(x)\leqq g_j(\bar x)$, $j\in A(\bar x)$, by Lemma \ref{lema3}, we have
$\scalpr{\nabla g_j(\bar x)}{x-\bar x}\leqq0$
 for all $j\in A(\bar x)$, which implies that the direction $x-\bar x$ is  critical.
Using the assumptions of the theorem, we obtain that  there exist vector multipliers $\lambda$ and $\mu$ with non-negative components such that 
\[
\mu_j g_j(\bar x)=0,\; j=1,...,m,\quad  \scalpr{\nabla L(\bar x)}{x-\bar x}=0,\quad  L^{\pr\pr}(\bar x,x-\bar x)\geqq 0.
\]
Therefore, $\mu_j=0$, when $j\notin A(\bar x)$. Using that the direction $x-\bar x$ is critical, we obtain
\[
\scalpr{\nabla L(\bar x)}{x-\bar x}=\sum_{i=1}^n\lambda_i\scalpr{\nabla f_i(\bar x)}{x-\bar x}+
\sum_{i\in A(\bar x)}\mu_i\,\nabla g_i(\bar x)(x-\bar x)\leqq 0.
\]
\noindent
Hence, 
\[
\begin{array}{l}
\lambda_i\scalpr{\nabla f_i(\bar x)}{x-\bar x}=0,\quad \forall \; i=1,2,\ldots,n \\
\mu_j\scalpr{\nabla g_j(\bar x)}{x-\bar x}=0,\quad \forall\; j\in A(\bar x).
\end{array}
\]
 Then, $\scalpr{\nabla f_i(\bar x)}{x-\bar x}=0$ for all indexes $i$ with $\lambda_i>0$, and $\scalpr{\nabla g_j(\bar x)}{x-\bar x}=0$ when $\mu_j>0$. It follows from second-order pseudoconvexity that %\begin{equation}\label{11}
$f_i^{\pr\pr}(\bar x,x-\bar x)<0$ for all $i=1,2,\ldots,n$ such that $\lambda_i\ne 0$. It follows from quasiconvexity of $g$ that 
\[
g_j^{\pr\pr}(\bar x,x-\bar x) = \lim_{t\to+0}\frac{g_j(\bar x+t(x-\bar x))-g_j(\bar x)-t \scalpr{\nabla g_j(\bar x)}{x-\bar x}}{t^2/2}\leqq 0.
\]
%\end{equation}
\noindent
 for all  $j\in A(\bar x)$  with $\mu_j>0$. There exist Lagrange multipliers $\lambda_i$ with strictly positive values, because $\lambda\ne 0$.
We conclude from here that
\[
L^{\pr\pr}(\bar x,x-\bar x)=\sum_{i=1}^n\lambda_i f_i^{\pr\pr}(\bar x,x-\bar x)
+\sum_{j\in A(\bar x),\;\mu_j>0}\mu_j\,g_j^{\pr\pr}(\bar x,x-\bar x)<0,
\]
%By quasiconvexity $f_i(\bar x+t(x-\bar x))\leqq f_i(\bar x)=0$
%for all $i\in A(\bar x)$ and for all sufficiently small $t\in[0,1]$.
%We conclude from here that  $L^{\pr\pr}(\bar x,x-\bar x)<0$
\noindent
which is a contradiction.
\end{proof}

The following example shows how Theorem \ref{th1} can be applied  in practical problems.

\begin{example}\label{ex1}
Consider the problem

\medskip
Minimize$\quad f(x)=(f_1(x),f_2(x))\quad$subject to $g(x)\leqq0$,

\medskip
\noindent
where the function $f_1:\R^2\to\R$ is defined as follows:
\[
f_1(x)=\left\{\begin{array}{rl}
x_1^2+x_2^2, & x_1\ge 0,\; x_2\ge 0, \\
x_2^2-x_1^2, & x_1\le 0,\; x_2\ge 0, \\ 
-x_1^2- x_2^2, & x_1\le 0,\; x_2\le 0, \\ 
x_1^2-x_2^2, & x_1\ge 0,\; x_2\le 0,
\end{array}\right.
\]
$f_2:\R^2\to\R$ and $g:\R^2\to\R$ are the functions of two variables 
\[
f_2(x)=-x_1-x_2+\sqrt{(x_1-x_2)^2+4},\quad g(x_1,x_2)=-x_1-x_2.
\]

The function $f_1$ is second-order pseudoconvex and continuously differentiable, $f_2$ is convex, and $g$ is linear. 
%To prove the second-order pseudoconvexity of $f_1$, we should consider  .Да изберем $\bar x=(0,0)$. Целевата функция е диференцируема в $\bar x$ и $\nabla f(\bar x)=(0,0)$. Също така тя е диференцируема от втори ред в $\bar x$ по всяко направление. Ще докажем, че е псевдоизпъкнала от втори ред в $\bar x$. Нека $f(y)<f(\bar x)$. Разглеждаме поотделно различните случаи $y$ да бъде точка от четирите квадранта. Например, нека $y$ е от втори квадрант. Следователно $y_2^2-y_1^2<0$. Оттук получяваме, че $\nabla f(\bar x)(y-\bar x)=0$, $f^{\pr\pr}(\bar x,y-\bar x)=2(y_2^2-y_1^2)<0$. Аналогично се разглежда четвърти квадрант и останалите случаи.

The points $x=(x_1,x_2)$ such that $x_1=x_2\ge 0$ are weakly effective. Indeed, let $x_1>x_2\geqq 0$, or $x_1+x_2>0$, $x_1>0$, $x_2<0$, or $x_1+x_2=0$, $x_1>0$. Then the direction $d=(d_1,d_2)$ such that $d_1+d_2=0$, $d_1<0$ ensure $\scalpr{\nabla f_1(x)}{d}<0$, $\scalpr{\nabla f_2(x)}{d}<0$. Let $x_2>x_1\geqq 0$, or $x_1+x_2>0$, $x_1<0$, $x_2>0$, or $x_1+x_2=0$, $x_2>0$. Then the direction $d=(d_1,d_2)$ such that $d_1+d_2=0$, $d_2<0$ ensure $\scalpr{\nabla f_1(x)}{d}<0$, $\scalpr{\nabla f_2(x)}{d}<0$. Therefore, these points are not weakly effective. Let $x_1=x_2\geqq 0$, $d_1+d_2\geqq 0$. Then $f_1(x+td)>f_1(x)$ for all sufficiently small positive numbers $t$. Let $x_1=x_2\geqq 0$, $d_1+d_2<0$. Then $f_2(x+td)>f_2(x)$ for all sufficiently small positive $t$. 
Suppose that $x_1=x_2=0$. If $d_1\geqq 0$, $d_2\geqq 0$, $d\ne 0$, or $d_1>0$, $d_2<0$, $d_1+d_2>0$, or $d_1<0$, $d_2>0$, $d_1+d_2>0$, then 
$f_1(x+td)>f_1(x)$ for all sufficiently small positive $t$. 
%Let $x_1=x_2=0$, $d_1>0$, $d_2<0$, $d_1+d_2>0$. Then $f_1(x+td)>f_1(x)$ for all sufficiently small positive $t$. Let $x_1=x_2=0$, $d_1<0$, $d_2>0$, $d_1+d_2>0$. Then $f_1(x+td)>f_1(x)$ for all sufficiently small positive $t$. 
If 
%$x_1=x_2=0$, 
$d_1+d_2=0$, then $f_1(x+td)=f_1(x)$, but $f_2(x+td)>f_2(x)$  for all sufficiently small positive numbers $t$. 

Let us find the weakly efficient points applying Theorem \ref{th1}. The points,  which satisfies the equations 
\[
\begin{array}{l}
\lambda_1\nabla f_1(x)+\lambda_2\nabla f_2(x)+\mu\nabla g(x)=0,\quad\mu g(x)=0,\\ 
 \lambda_1\geqq 0,\quad\lambda_2\geqq 0,\quad\mu\geqq 0,\quad (\lambda_1,\lambda_2)\ne(0,0)
\end{array}\]
 have the form $\bar x=(t,t)$, $t\geqq 0$ with multipliers $\lambda_1=1$, $\lambda_2=2t$, and $\mu=0$. The critical directions at $\bar x=(0,0)$ are $d=(d_1,d_2)$, where $d_1+d_2\geqq 0$. The critical directions at $\bar x=(t,t)$, $t>0$ are $d=(d_1,d_2)$, where $d_1+d_2=0$.
The second-order conditions are satisfied also. Then, it follows from the sufficient conditions from Theorem \ref{th1} that $\bar x$ is weak minimizer. 
\end{example}

This example cannot be solved with the first-order conditions, because $f_1$ is not pseudoconvex.

\begin{proposition}
If the active constraints satisfy Mangasarian-Fromovitz constraint qualification, then the second-order Mangasarian-Fromovitz constraint qualification holds.
\end{proposition}
\begin{proof}
%{\bf Proof.}
Suppose that  Mangasarian-Fromovitz constraint qualification is satisfied. Then, SOMFCQ holds with $K=\emptyset$.
\end{proof}

Recall that the closed convex hull of the tangent cone is called the pseudotangent cone. Let us consider the linearizing cone 
 \[
L(\bar x)=\{d\in\R^n\mid\scalpr{\nabla g_i(\bar x)}{d}\le 0,\; i\in A(\bar x)\}.
\] 
It is said that the Guignard constraint qualification \cite{gui69} is satisfied, iff the pseudotangent cone $P(S,x)$ of the feasible set $S$ at some feasible point $x$ coincides with the linearizing cone of the feasible set $L(x)$ at $x$.

The following example shows the possibility that  SOMFCQ is satisfied, but the Mangasa\-ri\-an-Fromovitz and Guignard constraint qualifications do not hold:
 
\begin{example}
Consider the following example:
 
\noindent
Minimize $f=x^4_1+x^4_2$\\
subject to the constraints $\quad g_1(x_1,x_2)\leqq 0,\quad g_2(x_1,x_2)\leqq0,$
where $g_1$ and $g_2$ are the functions
\[
g_1=\left\{
\begin{array}{rc}
x_1^2,\; & x_1\geqq 0 \\
-x_1^2, & x_1<0
\end{array}\right.\quad
\quad g_2=\left\{
\begin{array}{rc}
x_2^2, \; & x_2\geqq 0 \\
-x_2^2, & x_2<0
\end{array}\right.\quad
\]
\medskip
\noindent
The feasible set is $S=\{x=(x_1,x_2)\mid x_1\le 0,\; x_2\le 0\}$. The point $\bar x=(0,0)$ is a global solution and $A(\bar x)=\{1,2\}$. The pseudotangent cone $P(S,\bar x)$ coincides with the feasible set, and the linearizing cone $L(\bar x)$ coincides with the whole space $\R^2$. 
The Guignard CQ $L(\bar x)\subset P(S,\bar x)$ is not satisfied. Mangasarian-Fromovitz CQ is not satisfied also. On the other hand, the second-order Mangasarian-Fro\-mo\-vitz CQ is satisfied.
\end{example}

It is easy to verify Mangasarian-Fromovitz CQ. Another CQ, which is easy verified, is the Slater CQ. It is said that the Slater CQ is satisfied, if the constraint functions $g_j$, $j\in A(\bar x)$ are pseudoconvex, and there exists a point $x^0$ such that $g_j(x^0)<0$, $j\in A(\bar x)$.
We introduce another constraint qualification. 

\begin{definition}
We say that the constraint functions satisfy the second-order Slater CQ, iff the functions $g_j$, $j\in A(\bar x)$   are second-order pseudoconvex, and there exists a point $x^0$ such that $g_j(x^0)<0$, $j\in A(\bar x)$.
\end{definition}

It is easy to see that every problem which satisfies Slater CQ satisfies the second-order Slater CQ also.

It follows from second-order pseudoconvexity that,  if $\bar x$ is locally effective and $x^0$ satisfies second-order Slater CQ, then the direction $d=x^0-\bar x$ satisfies the SOMFCQ. Therefore, we can replace the SOMFCQ by the second-order Slater CQ in Theorem \ref{th1}.

\section{Lagrange Multiplier Characterizations of the Solution Set}
\label{s5}
In this section, we derive Lagrange multiplier characterizations of the solution set of a scalar problem with second-order pseudoconvex objective function and quasiconvex inequality constraints. We suppose that the problem has multiple solutions and one of them  $\bar x$, and the Lagrange multipliers, which satisfy KKT necessary conditions, are known.

Consider the problem with inequality constraints

\bigskip
Minimize $\quad f(x)\quad$ subject to $\quad x\in X,\quad g(x)\leqq 0, \; i=1,2,...,m$,\hfill  ${\rm (P)}$  
\bigskip

\noindent
where $f:\Gamma\to\R$ and $g:\Gamma\to\R^m$ are defined on some open set $\Gamma\subseteq\R^s$, $X$ is a convex subset of $\Gamma$, not necessarily open. Let  $S$ be the feasible set.

Suppose that $C\subseteq\R^s$ is a cone. Then, the cone 
\[
C^*:=\{c\in\R^s\; \mid \; \scalpr{c}{x}\le 0 \quad\textrm{for all}\quad x\in C \}
\]
 is said to be the negative polar cone of $C$.  Let $T_X(x)$ be the tangent cone of the set $X$ at the point $x$. Then, its negative polar cone is called the normal cone $N_X(x)$.

\begin{definition}\cite{ber03}
It is said that the constraint functions satisfy generalized Man\-ga\-sa\-ri\-an-Fromovitz constraint qualification (in short, GMFCQ) at the point $\bar x$, iff there is a direction $y\in (N_X(\bar x))^*$ such that $\scalpr{\nabla g_i(\bar x)}{y}<0$ for all $i\in A(\bar x)$.
\end{definition}
The following necessary optimality conditions of Karush-Kuhn-Tucker type (in short, KKT conditions) are consequence of Proposition 2.2.1, Definition 2.4.1 and Proposition 2.4.1 in \cite{ber03}:

\begin{KKT}\label{KKT}
Let $\bar x$ be a local minimizer of the problem {\rm (P)}. Suppose that $f$, $g_i$, $i\in A(\bar x)$ are Fr\'echet differentiable on 
$\Gamma\subseteq\R^n$ at $\bar x$, $g_i$, $i\notin A(\bar x)$ are continuous at $\bar x$, the set $X\subseteq\Gamma$ is convex. Suppose additionally that GMFCQ holds at $\bar x$.  Then, there exists a Lagrange multiplier 
\[
\mu\in\R^m,\quad\mu=(\mu_1,\mu_2,...,\mu_m),\quad \mu_i\ge 0,\; i=1,2,\dots,m
\]
 such that
\[
\scalpr{\nabla f(\bar x)+\sum_{i\in I(\bar x)}\mu_i\nabla g_i(\bar x)}{x-\bar x}\ge 0,\;\forall\, x\in X,\quad \mu_i g_i(\bar x)=0\;\; \forall\; i=1,2,...,m.
\]
\end{KKT}

Denote  by $\tilde A(\bar x,\mu)$ the following index set
\[
\tilde A(\bar x,\mu):=\{i\in\{1,2,...,m\}\; \mid \; g_i(\bar x)=0,\;\mu_i>0\}
\]
and the set
\[
X_1(\mu):=\{x\in X\; \mid \; g_i(x)=0\;\forall i\in\tilde A(\bar x,\mu),\;\; g_i(x)\le 0\;\forall i\in\{1,2,\dots,m\}\setminus\tilde A(\bar x,\mu)\}.
\]

In the proofs of the theorem in this section, we apply the following two results from the paper by Ivanov \cite{jota4}:

\begin{lemma}[\cite{jota4}, Lemma 3.7]\label{lema7}
Let the functions $f$ and $g$ be differentiable and quasiconvex, $\bar x\in\bar S$ be a solution. Suppose that the set $X$ is convex,  GMFCQ is satisfied at $\bar x$ and KKT optimality conditions are satisfied at $\bar x$ with a multiplier $\mu$.
Then, $\bar S\subseteq X_1(\mu)$.
\end{lemma}

\begin{lemma}[\cite{jota4}, Lemma 2.6]\label{lema6}
Let $\Gamma\subseteq\E$ be an open convex set, $S\subseteq\Gamma$ be a convex one. Suppose that $f:\Gamma\to\R$ is a continuously differentiable quasiconvex function. Then exactly one of the following alternatives holds:

I) $\nabla f(x)\ne 0$ for all $x\in\bar S$ and the normalized gradient $\nabla f(x)/\norm{\nabla f(x)}$ is constant over the solution set $\bar S$;

II) $\nabla f(x)=0$ for all $x\in\bar S$.
\end{lemma}

Let $\Gamma\subseteq\E$ be an open set.
Recall that the following directional derivative of a Fr\'echet differentiable function $\f$  at a point $x\in\Gamma$ in direction $d\in\E$
\[
\sld f(x,d):=\liminf_{t\to 0^+}\,2t^{-2}\,[f(x+td)-f(x)-t\scalpr{\nabla f(x)}{d}]
\]
is usually called the second-order lower Dini directional derivative (or Peano derivative). Every Fr\'echet differentiable function has a second-order lower Dini derivative, eventually infinite. 

Consider the following sets:
\begin{align}
S_1:=\{x\in X_1(\mu)\; \mid\;\scalpr{\nabla f(\bar x)}{x-\bar x}=0,\;\exists p(x)>0:\nabla f(x)=p(x)\nabla f(\bar x),\nonumber\\
\exists f^{\pr\pr}(x,\bar x-x),\;f^{\pr\pr}(x,\bar x-x)=0\};\nonumber
\end{align}
\begin{align}
S_2:=\{x\in X_1(\mu)\; \mid\;\scalpr{\nabla f(\bar x)}{x-\bar x}\leqq 0,\;\exists p(x)>0:\nabla f(x)=p(x)\nabla f(\bar x),\nonumber\\
\exists f^{\pr\pr}(x,\bar x-x),\;f^{\pr\pr}(x,\bar x-x)\geqq 0\};\nonumber
\end{align}
\begin{align}
S_3:=\{x\in X_1(\mu)\; \mid\;\scalpr{\nabla f(\bar x)}{x-\bar x}=0,\;\exists p(x)>0:\nabla f(x)=p(x)\nabla f(\bar x),\nonumber\\
\exists f^{\pr\pr}(x,\bar x-x),\;f^{\pr\pr}(x,\bar x-x)=\sld f(\bar x,x-\bar x)\};\nonumber
\end{align}
\begin{align}
S_4:=\{x\in X_1(\mu)\; \mid\;\scalpr{\nabla f(\bar x)}{x-\bar x}\leqq 0,\;\exists p(x)>0:\nabla f(x)=p(x)\nabla f(\bar x),\nonumber\\
\exists f^{\pr\pr}(x,\bar x-x),\;f^{\pr\pr}(x,\bar x-x)\geqq\sld f(\bar x,x-\bar x)\};\nonumber
\end{align}
\begin{align}
S_5:=\{x\in X_1(\mu)\; \mid\;\scalpr{\nabla f(\bar x)}{x-\bar x}=0,\;\exists p(x)>0:\nabla f(x)=p(x)\nabla f(\bar x),\nonumber\\
\exists f^{\pr\pr}(x,\bar x-x),\;f^{\pr\pr}(x,\bar x-x)=\sld f(\bar x,x-\bar x)=0\}. \nonumber
\end{align}
We do not suppose that the existence of the second-order directional derivatives is gua\-ranteed. This case is more general than the case, which includes the assumption that $f$ is second-order differentiable in every direction. For example, if $f^{\pr\pr}(x,\bar x-x)$ does not exist, then $x\notin \tilde A$.

\begin{theorem}\label{th3}
Let the function $f$ be continuously differentiable and second-order pseudoconvex, $g$ be differentiable and quasiconvex, $X$ be a convex set
 $\bar x\in\bar S$ be a known solution of {\rm (P)}, $\Gamma\subseteq\E$ be an open convex set, and GMFCQ be satisfied. Suppose that $\mu\in\R^m_+$ is a known vector Lagrange multiplier, which fulfills KKT conditions. Then,
%If the function $f$ is continuously differentiable and second-order pseudoconvex on $\Gamma$, then
\[
\bar  S=S_1=S_2=S_3=S_4=S_5.
\]
\end{theorem}
\begin{proof}
It is obvious that $S_5\subseteq S_1\subseteq S_2$ and $S_5\subseteq S_3\subseteq S_4$.

Consider the case, when $\nabla f(\bar x)\ne 0$.
We prove that  $\bar  S\subseteq S_5$. Suppose that $x\in\bar S$. Therefore, $f(x)=f(\bar x)$. By Lemma \ref{lema2} the function $f$ is quasiconvex. By quasiconvexity, the level sets of $f$ are convex. Therefore, $\bar S$ is also convex. Thus, we obtain that
\[
f[\bar x+t(x-\bar x)]=f[x+t(\bar x-x)]=f(\bar x)\quad\textrm{for all}\quad t\in [0,1].
\]
It follows from here that $\scalpr{\nabla f(\bar x)}{x-\bar x}=0$. We can prove using similar arguments,  interchanging $\bar x$ and $x$, that $\scalpr{\nabla f(x)}{\bar x-x}=0$. We conclude from the definition of the second-order derivative that $f^{\pr\pr}(x,\bar x-x)$ exists  and
\[
f^{\pr\pr}(x,\bar x-x)=\lim_{t\to +0}\frac{f[x+t(\bar x-x)]-f(x)-t\scalpr{\nabla f(x)}{\bar x-x}}{t^2/2}=0.
\]
We can prove using similar arguments,  interchanging $\bar x$ and $x$, that $\sld f(\bar x,x-\bar x)=0$.
It follows from Lemma \ref{lema7} that $\bar S\subseteq X_1(\mu)$. By Lemma \ref{lema6} and the assumption $\nabla f(\bar x)\ne 0$, we conclude that $\nabla f(x)\ne 0$ and $\nabla f(x)/\norm{\nabla f(x)}=\nabla f(\bar x)/\norm{\nabla f(\bar x)}$. Then, we obtain that $\nabla f(x)=p(x)\nabla f(\bar x)$, where $p(x)=\norm{\nabla f(x)}/\norm{\nabla f(\bar x)}$. It is obvious that $p(x)>0$. It follows from all these arguments that $x\in S_5$.

We prove that $S_4\subseteq S_2$. Let $x\in S_4$. Therefore, there exists the second-order derivative $f^{\pr\pr}(x,\bar x-x)$ and 
\begin{equation}\label{24}
f^{\pr\pr}(x,\bar x-x)\geqq\sld f(\bar x,x-\bar x).
\end{equation}
 Since the set $X$ is convex, and the functions $g_i$ are quasiconvex, then the feasible set $S$ is convex. It follows from here, by $x\in S$ and $\bar x\in\bar S$, that $f[\bar x+t(x-\bar x)]\geqq f(\bar x)$ for all $t\in [0,1]$. Using that $\scalpr{\nabla f(\bar x)}{x-\bar x}\leqq 0$, we infer from the definition of the second-order lower derivative that
\[
\sld f(\bar x,x-\bar x)=\liminf_{t\to +0}\frac{f[\bar x+t(x-\bar x)]-f(\bar x)-t\scalpr{\nabla f(\bar x)}{x-\bar x}}{t^2/2}\geqq 0.
\]
Then, (\ref{24}) implies that $f^{\pr\pr}(x,\bar x-x)\geqq 0$. Therefore $x\in S_2$.

At last, we prove that $S_2\subseteq\bar S$. Let $x\in S_2$. Assume the contrary that $x\notin\bar S$. Hence $f(\bar x)<f(x)$. By
second-order pseudoconvexity, we obtain that 
\[
\scalpr{\nabla f(x)}{\bar x-x}\le 0.
\]
 Then, the condition $\nabla f(x)=p\nabla f(\bar x)$, $p>0$ leads us to the conclusion that 
\[
\scalpr{\nabla f(\bar x)}{x-\bar x}\geqq 0.
\]
 The last inequality together with $\scalpr{\nabla f(\bar x)}{x-\bar x}\leqq 0$ implies that $\scalpr{\nabla f(\bar x)}{x-\bar x}=0$. Again from the equation $\nabla f(x)=p\nabla f(\bar x)$, $p>0$, we get that 
 $\scalpr{\nabla f(x)}{\bar x-x}=0$.  Then, by second-order pseudoconvexity of $f$ and $f(\bar x)<f(x)$, we conclude that $f^{\pr\pr}(x,\bar x-x)<0$, which contradicts the relation $x\in S_2$.

The proof in case, when $\nabla f(\bar x)=0$, can be obtained using the same scheme. It is simpler, than the presented proof.
\end{proof}

\begin{example}
Consider the function $f:\R^2\to\R$ and the problem

\medskip
Minimize $f(x)$ subject to $x\in X$, $g(x)=g(x_1,x_2)\leqq 0$,
\medskip

\noindent
where $X=\{x=(x_1,x_2)\mid x_1\geqq 0,\;-\infty<x_2<+\infty\}$,  $g(x_1,x_2)=-x_1$, and 
\[
f(x)=\left\{\begin{array}{rl}
x_1^2+x_2^2, & x_1\ge 0,\; x_2\ge 0, \\
x_2^2, & x_1\le 0,\; x_2\ge 0, \\ 
-x_2^2, & x_1\le x_2\le 0, \\ 
-x_1^2, & x_2\le x_1\le 0, \\
x_1^2, & x_1\ge 0,\; x_2\le 0.
\end{array}\right.
\]
The solution set is $\bar S:=\{x\in\R^2\, \mid \, x_1=0,\; x_2\leqq 0\}$. The constraint function $g$ is linear; therefore quasiconvex and differentiable. GMFCQ is satisfied, because 
\[
(N(X))^*=X^*=\{y=(y_1,y_2)\mid y_1\le 0,\; y_2=0\}.
\]
The objective function is not differentiable only over the set 
\[
B:=\{(x_1,x_2)\mid x_1=x_2<0\}.
\]
 In particular, it is differentiable over the feasible set $S$. It is also continuously differentiable and second-order pseudoconvex over $S$. The function $f$ is not pseudoconvex. We can apply Theorem \ref{th3}, because Lemma 2.6 in Ref. \cite{jota4} is satisfied, if the objective function is continuously differentiable over the feasible set. Let us take $\bar x=(0,0)$. The Lagrange multiplier, which satisfies KKT conditions at $\bar x$ is $\mu=0$. Only the points from the feasible set such that  $x_1=0$, $x_2\leqq 0$ are solutions of the equation $\nabla f(x)=p\nabla f(\bar x)=(0,0)$. It is easy to check that $\scalpr{\nabla f(\bar x)}{x-\bar x}=0$ and $f^{\pr\pr}(\bar x,x-\bar x)=0$ hold if $x_1=0$, $x_2\leqq 0$. Then, it follows from Theorem \ref{th3} that $\bar S:=\{x\in\R^2\, \mid \, x_1=0,\; x_2\leqq 0\}$. 
\end{example}

\section{Necessary and sufficient conditions for an efficient solution}
\label{s4}

In this section, we derive necessary and sufficient optimality conditions for  efficiency in the vector problem {\rm{(VP)}}.

We introduce the following definition:

\begin{definition}
We call the problem {\rm (VP)} second-order KKT-pseudoconvex at the point $x\in S$, iff the following implication holds:
\[
\left.
\begin{array}{l}
y\in S, \\
f(y)\le f(x)
\end{array}\right]
\;\Rightarrow\;\left[
\begin{array}{l}
\scalpr{\nabla f_i(x)}{y-x}\leqq 0,\; i=1,2,\dots,n \\
\scalpr{\nabla f_i(x)}{y-x}= 0\quad\textrm{implies}\quad f^{\pr\pr}_i(x;y-x)<0,\\
\scalpr{\nabla g_j(x)}{y-x}\leqq 0\;\textrm{ for all }\; j\in A(x), \\
\scalpr{\nabla g_j(x)}{y-x}=0,\; j\in A(x)\quad\textrm{imply}\quad g_j^{\pr\pr}(x;y-x)\leqq 0,
\\
\end{array}
\right.
\]
provided that all necessary derivatives exist. We call the problem {\rm (VP)} second-order  KKT-pseudoconvex, iff it is second-order  KKT-pseudoconvex at each $x\in S$.
\end{definition}

In the case when $n=1$ this notion reduces to the notion of KT pseudoconvex scalar problem (see Ivanov \cite{jogo1}). %\cite{jogo}).

\begin{theorem}\label{th2}
Let the point $\bar x$ be feasible for the problem {\rm (VP)}. Suppose that the set $X$ is convex, the vector function $f$ and the functions $g_j$, $j\in A(\bar x)$ are continuously differentiable, the functions $g_j$, $j\notin A(\bar x)$ are continuous at $\bar x$, $f$ and $g$ are second-order directionally differentiable at $\bar x$ in every critical direction $d$, 
the problem (VP) is second-order KKT-pseudoconvex. 

Let $\bar x$ be a Pareto minimizer. Then, for every direction $d$ there exist Lagrange multipliers $\lambda$, $\mu$ with non-negative components such that $(\lambda,\mu)\ne (0,0)$ and conditions (\ref{9}), if $d$ is not critical,  (\ref{10}), if $d$ is critical, are satisfied. Suppose additionally that SOMFCQ holds. Then, then for every direction $d$, satisfying SOMFCQ, there exist 
vector multipliers $\lambda$, $\mu$ with non-negative components such that (\ref{9}) and (\ref{10}) hold with $\lambda\ne 0$. 

Conversely, let for every critical direction $d$ there exist Lagrange multipliers $\lambda$, $\mu$ with non-negative components such that $\lambda\ne 0$ and conditions (\ref{10}) are satisfied.
Then, $\bar x$ is a globally efficient solution.
\end{theorem}

\begin{proof}
The necessity follows from Theorem \ref{th1} and the fact that every efficient point is weakly efficient.

We prove the sufficiency. Suppose that $x\in S$ is a second-order Karush-Kuhn-Tucker stationary point, but it is not a global Pareto minimizer. Therefore, there exists $y\in S$ such that $f(y)\le f(x)$. It follows from second-order KKT-pseudoconvexity of (VP) that 
\begin{equation}\label{11}
\begin{array}{l}
\scalpr{\nabla f_i(x)}{y-x}\leqq 0,\quad \forall i=1,2,\dots, n, \\
\scalpr{\nabla g_j(x)}{y-x}\leqq 0,\quad \forall j\in A(x).
\end{array}
\end{equation}
 We obtain from here that the direction $y-x$ is critical at the point $x$. According to second-order Karush-Kuhn-Tucker conditions we have $\scalpr{\nabla L(x)}{y-x}=0$. Then, it follows from (\ref{11}) that 
\[
\lambda_i\scalpr{\nabla f_i(x)}{y-x}=0,\quad \mu_j\scalpr{\nabla g_j(x)}{y-x}=0,\quad\forall i=1,2,\dots,n,\;\forall j\in A(\bar x).
\]
 We conclude from here that
\[
\scalpr{\nabla f_i(x)}{y-x}=0, \quad \scalpr{\nabla g_j(x)}{y-x}=0
\]
for all $i$ and $j$ such that $\lambda_i>0$ and $\mu_j>0$. Then it follows from second-order KKT-pseudoconvexity of (VP) that 
\[
f^{\pr\pr}_i(x;y-x)<0,\quad g_j^{\pr\pr}(x;y-x)\leqq 0
\]
 for all $i$ and $j\in A(x)$ such that $\lambda_i>0$ and $\mu_j>0$.
There exist Lagrange multipliers $\lambda_i$ with strictly positive values, because $\lambda\ne 0$.
We conclude from here that
%By $\lambda\ne 0$ we get the inequality 
\[
L^{\pr\pr}(\bar x,y-x)=\sum_{i=1}^n\lambda_i f_i^{\pr\pr}(\bar x,y-x)+ \sum_{j\in A(\bar x)}\mu_j g_j^{\pr\pr}(\bar x,y-x)< 0,
\]
%\noindent
which contradicts the second-order Karush-Kuhn-Tucker condition (\ref{10}).
\end{proof}

%The problem from Example \ref{ex1} is really KT-pseudoconvex. Therefore $\bar x=(0,0)$ is a Pareto minimizer.

\end{document}